\documentclass[10pt,twoside]{amsart}

\usepackage{graphicx}
\usepackage{color}
\usepackage{amssymb}
\usepackage{latexsym}
\usepackage[T1]{fontenc}
\usepackage[utf8]{inputenc}
\usepackage[english]{babel}
\usepackage{amsmath}
\usepackage{amsthm}
\usepackage{amsfonts}
\usepackage{amsxtra}
\usepackage{enumerate}
\usepackage{verbatim}
\usepackage{color}

\makeatletter

\@addtoreset{equation}{section} \makeatother
\newtheorem{theorem}{Theorem}[section]
\newtheorem{remark}[theorem]{Remark}
\newtheorem{lemma}[theorem]{Lemma}
\newtheorem{proposition}[theorem]{Proposition}
\newtheorem{corollary}[theorem]{Corollary}
\newtheorem{definition}[theorem]{Definition}
\newtheorem{example}[theorem]{Example}

\def\1{\mathbf{1}}
\def\:{\lrcorner}
\def\#{\sharp}

\def\l{\lambda}

\def\g{\gamma}

\def\qed{\ensuremath{\quad\Box\quad}}

\def\inv#1{\raise.1em\hbox to 0pt{$^{-1}$\hss}_{#1}\;}

\def\v{\noindent}

\newcommand{\bean}{\begin{eqnarray*}}
\newcommand{\eean}{\end{eqnarray*}}
\newcommand{\benu}{\begin{enumerate}}
\newcommand{\eenu}{\end{enumerate}}
\newcommand{\eea}{\end{eqnarray}}
\newcommand{\bea}{\begin{eqnarray}}

\newtheorem{Theorem}{Theorem}

\newtheorem{Definition}{Definition}

\newcommand{\be}{\begin{equation}}
\newcommand{\ee}{\end{equation}}

\newcommand{\R}{{\mathbb R}}

\newcommand{\ben}{\begin{enumerate}}
\newcommand{\een}{\end{enumerate}}
\newcommand{\bit}{\begin{itemize}}
\newcommand{\eit}{\end{itemize}}
\newcommand{\edoc}{\end{document}}

\newcommand{\bdefi}{\begin{definition}}
\newcommand{\btheo}{\begin{theorem}}
\newcommand{\bprop}{\begin{proposition}}
\newcommand{\brema}{\begin{remark}}
\newcommand{\bcoro}{\begin{corollary}}
\newcommand{\blemm}{\begin{lemma}}
\newcommand{\bexam}{\begin{example}}

\newcommand{\edefi}{\end{definition}}
\newcommand{\etheo}{\end{theorem}}
\newcommand{\eprop}{\end{proposition}}
\newcommand{\erema}{\end{remark}}
\newcommand{\ecoro}{\end{corollary}}
\newcommand{\elemm}{\end{lemma}}
\newcommand{\eexam}{\end{example}}

\begin{document}

\title{Asymptotic flexibility of globally hyperbolic manifolds}

\author[]{Olaf M\"uller}
\address{Fakult\"at f\"ur Mathematik, Universit\"at Regensburg}
\email{olaf.mueller@mathematik.uni-regensburg.de}

\parindent=5mm
\date{\today}

\begin{abstract}
In this short note, a question of patching together globally hyperbolic manifolds is adressed which appeared in the context of the construction of Hadamard states.
\end{abstract}

\maketitle


\v Often, for a normally hyperbolic field theory (as Maxwell or Klein-Gordon theory) on a globally hyperbolic manifold, one wishes to construct Hadamard states. Those are complex-linear functionals on the Weyl algebra (which, in turn, is a certain subalgebra on the algebra of smooth complex functions on a space of solutions to some field equation) satisfying additional properties, for details see \cite{mR}. The crucial point for this note consists solely in the fact that, while the Hadamard property can be defined locally, to every state defined in, say, an open causal and thus globally hyperbolic neighborhood of some Cauchy surface we can associate a Hadamard state in all of the manifold. This procedure is called {\em propagation of the state} (to the future or the past). In ultrastatic globally hyperbolic manifolds, there is an easy and very explicit method for the construction of Hadamard states. Now if we know that we can modify the past of a Cauchy surface of a given manifold $(M,g)$ in a way that the modified metric $(M, \tilde{g})$ is asymptotically ultrastatic (while staying globally hyperbolic) then we can define a Hadamard state in the past in $\tilde{g}$ and propagate it to the future. According to what has been said above, it will stay Hadamard. Then we consider the Hadamard state in the future in which $g$ and $\tilde{g}$ coincide and propagate it back to the past of the original metric $g$. The state we have constructed this way is Hadamard for the original metric. Now the question arises if this construction can be performed for every globally hyperbolic manifold. Sometimes a slightly different procedure is done in which, for a given globally hyperbolic manifold $(M,g)$, another one $(M,\tilde{g})$ is constructed which is ultrastatic in the past and contains an open neighborhood $N$ of a Cauchy surface $S$ of $(M,g)$ (cf. \cite{rV}, for example). But the size of $N$ cannot be controlled due to the proof which works by Fermi coordinates around $S$. Thus the construction, although very useful for showing the existence of Hadamard spaces, leaves questions from Lorentzian geometry involving concepts as geodesic completeness unanswered. The following result answers the above question in the affirmative. 

\begin{Definition}
Two globally hyperbolic manifolds $(M,g)$ and $(N,h)$ are called {\bf future-isometric} (resp. {\bf past-isometric}) iff there is a Cauchy hypersurface $S$ of $(M,g)$ and $T$ of $(N,h)$ such that $I^+(S)$ is isometric to $I^+ (T)$ (resp. $I^-(S)$ is isometric to $I^-(T)$). Let $J(g,h)$ be the set of globally hyperbolic manifolds past-isometric to $g$ and future-isometric to $h$. Any metric in $J(g,h)$ is called an {\bf asymptotic join of $g$ and $h$}.   
\end{Definition}

\v We define a binary symmetric relation $P$ of past-isometry (resp. $F$ of future-isometry) which is moreover transitive, as for two different Cauchy surfaces one can find a third one in the past (resp. future) of both. We will first prove another result used in the proof of the second theorem: 


\begin{Theorem}
\label{erstes}
Let $\l \in C^{\infty} ( \R \times N) $  and let $g:= - dt^2 + g_t $ be a Lorentzian metric on $\R \times N$, where each $g_t$ is a Riemannian metric on $\{t \} \times N$. Then there is an $ f \in C^{\infty} ( \R \times N) $ such that $(M:= \R \times N , h:=  - \l dt^2 + f g_t ) $ is globally hyperbolic. If $\l_s = \l_u ,g_s = g_u$ for any two $ s,u \in (- \infty, 0)$, and if $ (a, \infty) \times N  $ is already globally hyperbolic, then $f$ can be chosen such that $f_s = f_u$ for any two $s, u \in (- \infty,0)$ as well and equal to one on $ (a, \infty) $. 
\end{Theorem}

\v {\bf Proof.} First choose a smooth function $j$ on $S := t^{-1} (\{ 0 \})$ such that $j g_0$ is complete. Now we want to have $I^{\pm}_{h} (t,x) \cap S \subset B^S_t(x) $ (which ensures global hyperbolicity). As we can parametrize any causal curve $c$ as $ c(t) = (t, k(t) ) $ and as for the resulting curve $k$ holds $f_t g_t (\dot{k}, \dot{k} ) \leq \l_t$, it is sufficient that $ f_t g_t \geq {\rm max} \{ 1, \l \} \cdot j g_0 $. By compactness of the Euclidean sphere, there is a continuous function $\underline{f}$ satisfying this inequality, so we can choose a smooth function $f \geq \underline{f}$ with this property as well. The additional property is now obvious as the choice of $\underline{f}$ was pointwise.   \hfill \qed 


\begin{Theorem}
Let $(M,g)$ and $(M,h) $be globally hyperbolic, let the Cauchy hypersurfaces of $ g  $ be diffeomorphic to those of $h$. Then $J(g,h)$ is nonempty. In particular, for any $(M,g)$ globally hyperbolic, there is a globally hyperbolic ultrastatic metric $u$ on $M$ such that $J(g,u)$ is nonempty.
\end{Theorem}

\bigskip

\v {\bf Proof.} Choose a metric splitting $(M,g) = (\R \times N, - s \cdot dt^2 + g_t) $ by a smooth Cauchy time function $t$ as in \cite{mo} and put $N := t^{-1} (1)$ and $S:= t^{-1} (0)$. Then choose a smooth positive functions $f$ on $M = \R \times N$ such that $ f \vert_{I^+ (S)} = 1$ and with $ f = s^{-1}$  on $I^-(N)$. Then, via $t$, the metric $g^{(1)} := f \cdot g$ splits as $(M, g^{(1)} ) = (-l dt^2 , f g_t )$ and $l= -1$ in $I^-(N)$. Moreover, $F(g,fg)$. Now, for a smooth monotonously increasing function $\psi : \R \rightarrow [0, \infty)$ with $\psi(r) = 0 \qquad \forall r \leq 0$, $\psi(r) = r \qquad \forall r \geq 1$, define a smooth function $\l$ and a Lorentzian metric $ k := - dt^2 + k_t$ as in Theorem \ref{erstes} by $\l_t := l_{\theta (t)}$ and $k_t := f_{\theta (t)} \cdot g_{\theta (t)}$. Note that $ \l_t $ and $k_t$ are constant for $t<0$. Then apply the first theorem to $(\l, k)$ and get a smooth functions $\phi $ on $M$ such that $(\R \times N, \g := - \l_t dt^2 + \phi_t k_t)$ is globally hyperbolic, and $\phi_t$ can be chosen equal to $1$ on $[1, \infty)$ and constant on $( - \infty, 0]$. Then $ F(g, \g) $, and $ P(\g, u)$ where $u$ is the ultrastatic metric $-dt^2 + f g_{-1}$. Therefore $\g \in J(g,u)$. If we have two different globally hyperbolic metrics $g$ and $h$ we construct the ultrastatic metrics $u_g = - dt^2 + k_0$ and $u_h = -dt^2 + k_1$ as above and intermediate between them via the metric $u_{gh} := -dt^2 + k_{\theta (t) }$ where $ k_r := r k_1 + (1-r) k_0  $ and $\theta : \R \rightarrow [0,1]$ smooth and monotonously increasing with $\theta ((- \infty, 0]) = \{ 0 \}$, $\theta ' (r) \neq 0 $ for all $ r \in (0,1)$, and $\theta ([1, \infty) = \{ 1 \}$. Then $u_{gh}$ is globally hyperbolic because it is stably causal and the causal diamonds $D(p,q)$ are compact as they are compact in every subset of an open covering defined by $A:= t^{-1} ((- \infty, 2/3))$ and $B:= t^{-1} ((1/3, \infty))$: If $t(q)  < 2/3$ this follows from comparison with the complete metric $(1 - \theta (2/3)) g_0$, in the other case $t(p) > 1/3$ it follows from the comparison with the complete metric $\theta (1/3) g_1$. Thus, as there is an open covering by two g.h. manifolds with a joint Cauchy surface, the entire manifold is g.h.   \hfill \qed 




\end{document}